\documentclass[12pt]{smfart}
\usepackage{smfthm}
\usepackage{smfenum}
\usepackage{amssymb}
\usepackage{amscd}
\usepackage[french]{babel}
\newcommand{\on}{\operatorname}

\newcommand{\Qp}{\mathbf{Q}_p}
\newcommand{\Zp}{\mathbf{Z}_p}

\newcommand{\ZZ}{\mathbf{Z}}
\newcommand{\NN}{\mathbf{N}}
\newcommand{\RR}{\mathbf{R}}

\newcommand{\bnrig}[2]{\mathbf{B}^{\dagger #1}_{\mathrm{rig} #2}}

\newcommand{\bdR}{\mathbf{B}_{\mathrm{dR}}}

\newcommand{\bdag}[1]{\mathbf{B}^{\dagger #1}}

\newcommand{\dnrig}[1]{\mathbf{D}^{\dagger #1}_{\mathrm{rig}}}

\newcommand{\dst}{\mathbf{D}_{\mathrm{st}}}
\newcommand{\dcris}{\mathbf{D}_{\mathrm{cris}}}
\newcommand{\ddR}{\mathbf{D}_{\mathrm{dR}}}

\renewcommand{\ddag}[1]{\mathbf{D}^{\dagger #1}}

\newcommand{\ndr}{\mathbf{N}_{\mathrm{dR}}}

\newcommand{\eps}{\varepsilon}
\newcommand{\ra}{\rightarrow}

\renewcommand{\phi}{\varphi}
\renewcommand{\projlim}{\varprojlim}

\newcommand{\Fil}{\mathrm{Fil}}

\renewcommand{\geq}{\geqslant}
\renewcommand{\leq}{\leqslant} 

\newcommand{\Iw}{\mathrm{Iw}}
\newcommand{\ddagrig}{D^\dag_{\mathrm{rig}}}

\title{Repr\'esentations de de Rham et normes universelles}
\alttitle{De Rham representations and universal norms}

\author{Laurent Berger}

\address{Harvard Dept of Mathematics  \\
      One Oxford Street \\
      Cambridge, MA 02138-2901 \\ USA}
\email{laurent@math.harvard.edu}
\urladdr{www.math.harvard.edu/\~{}laurent}


\date{11 Septembre 2003}
\begin{document}

\begin{abstract}
On calcule le module des normes universelles pour une repr\'esentation
$p$-adique de de Rham. Le calcul utilise la th\'eorie des 
$(\varphi,\Gamma)$-modules 
(la formule de r\'eciprocit\'e de Cherbonnier-Colmez)
et l'\'equation diff\'erentielle associ\'ee
\`a une repr\'esentation de de Rham.
\end{abstract}

\begin{altabstract}
We compute the module of universal norms for a  
de Rham $p$-adic representation. The computation uses the theory of 
$(\varphi,\Gamma)$-modules (Cherbonnier-Colmez's reciprocity formula) 
and the differential equation attached to a de Rham representation.
\end{altabstract}

\maketitle
\setcounter{tocdepth}{2}
\tableofcontents

\setlength{\baselineskip}{18pt}

\section*{Introduction}
Dans tout cet article, $p$ est un nombre premier,
$k$ est un corps parfait de caract\'eristique
$p$, et $K$ est une extension finie
du corps des fractions $F$ de l'anneau des vecteurs de Witt
sur $k$. On se fixe une cl\^oture alg\'ebrique $\overline{K}$ de $K$,
et on pose $G_K=\on{Gal}(\overline{K}/K)$. Rappelons qu'en utilisant
l'anneau de p\'eriodes $p$-adiques
$\bdR$, Fontaine a d\'efini la notion de
repr\'esentation de \emph{de Rham} de $G_K$. 
Si $V$ est une telle
repr\'esentation et si $L$ est une extension finie de $K$, on d\'efinit
$H^1_g(L,V)$ comme \'etant
l'ensemble des classes de cohomologie qui d\'eterminent
une extension 
$0 \ra V \ra E \ra \Qp \ra 0$ de repr\'esentations de $G_L$ telle que 
$E$ est une repr\'esentation de de Rham de $G_L$.

On \'ecrira $\mu_{p^n} \subset \overline{K}$ pour d\'esigner 
l'ensemble des racines $p^n$-\`emes de l'unit\'e, et on d\'efinit
$K_n=K(\mu_{p^n})$ ainsi que $K_{\infty} = \cup_{n=0}^{+\infty}
K_n$. On pose $H_K=\on{Gal}(\overline{K}/K_\infty)$ et  
$\Gamma_K=\on{Gal}(K_{\infty}/K)$. L'alg\`ebre d'Iwasawa est
l'alg\`ebre de groupe compl\'et\'ee $\Lambda_K=\Zp[[\Gamma_K]]$.

La cohomologie d'Iwasawa de $V$ est d\'efinie par $H^1_{\Iw}(K,V) = \Qp
\otimes_{\Zp} H^1_{\Iw}(K,T)$ o\`u $T$ est un r\'eseau de $V$ stable
par $G_K$ et $H^1_{\Iw}(K,T) = \projlim_n H^1(K_n,T)$ 
est la limite projective pour les applications 
{\og corestriction \fg} ce qui fait de $H^1_{\Iw}(K,V)$ un 
$\Qp \otimes_{\Zp} \Lambda_K$-module.
Par construction, on a des applications
$\on{pr}_{K_n,V}: H^1_{\Iw}(K,V) 
\ra H^1(K_n,V)$ et l'objet de cet
article est l'\'etude du module $H^1_{\Iw}(K,V)_g$ des {\og normes
  universelles \fg}, 
l'ensemble des $y \in H^1_{\Iw}(K,V)$ tels que pour tout $n \geq 0$ on
ait $\on{pr}_{K_n,V}(y) \in H^1_g(K_n,V)$.

Si $\Fil^1 V$ est la plus grande sous-repr\'esentation de $V$ 
dont tous les poids de Hodge-Tate sont $\geq 1$, alors  
$\Fil^1 V$ est de de Rham et 
on peut montrer que 
$H^1_{\Iw}(K,\Fil^1 V)_g = H^1_{\Iw}(K,\Fil^1 V)$.
Le r\'esultat principal de cet article est le suivant:
\begin{theo}\label{main_theo}
Si $V$ est une repr\'esentation de de Rham, alors
\begin{enumerate}
%\item $y \in  H^1_{\Iw}(K,V)_g$ si et seulement si pour tout $n \geq
%0$, on a $\on{pr}_{K_n,V}(y) \in H^1(K_n,\Fil^1 V)$;
\item si $V$ n'a pas de sous-quotient
stable par $H_K$, alors 
$H^1_{\Iw}(K,V)_g = H^1_{\Iw}(K,\Fil^1 V)$;
\item en g\'en\'eral, 
$H^1_{\Iw}(K,\Fil^1 V) \subset H^1_{\Iw}(K,V)_g$ et le quotient est un
$\Qp \otimes_{\Zp} \Lambda_K$-module de torsion.
\end{enumerate}
\end{theo}

La d\'emonstration est tr\`es similaire \`a celle qu'en a donn\'ee
Perrin-Riou dans \cite{BP00} pour les repr\'esentations cristallines,
et dans \cite{BP01} pour les repr\'esentations semi-stables,
du groupe de Galois d'un corps non-ramifi\'e
\footnote{au moins si les pentes 
de Frobenius sur $\dcris(V)$ 
sont enti\`eres, voir le d\'ebut des paragraphes
\cite[\S 3.2, \S 4.1]{BP00}} mais au lieu d'utiliser
son {\og Exponentielle \'elargie \fg}, on utilise les constructions de 
\cite{LB02} (qui redonnent l'exponentielle \'elargie d'ailleurs,
comme on le montre dans \cite{LB03}) et au lieu d'utiliser des
consid\'erations de {\og cran de la filtration \fg} et d'ordre, on
utilise la th\'eorie des $(\phi,\Gamma)$-modules qui encodent toutes ces
informations. Cela simplifie la d\'emonstration de Perrin-Riou et nous
permet de plus de l'\'etendre au cas des repr\'esentations de de Rham
du groupe de Galois d'un corps \'eventuellement ramifi\'e. Remarquons
que l'on n'utilise pas le fait que les repr\'esentations de de Rham
sont potentiellement semi-stables.

Indiquons bri\`evement d'o\`u provient cette conjecture. On renvoie \`a
l'article \cite{BP00} de Perrin-Riou pour plus d'informations. Si $E$
est une courbe elliptique d\'efinie sur $K$, on s'int\'eresse au
module des {\og normes universelles \fg} $\mathcal{N}_{K_\infty/K}(E) 
= \cap_{n = 0}^{+\infty} \on{Tr}_{K_n/K} E(K_n)$. 
Le calcul de ce module a \'et\'e fait 
(pour une courbe elliptique, une vari\'et\'e ab\'elienne 
ou m\^eme un groupe formel, et pour une $\Zp$-extension quelconque) dans
certains cas par Mazur \cite{BM72}, puis par Hazewinkel
\cite{MH74a,MH74b,MH77}, Schneider \cite{PS87} et Perrin-Riou
\cite{BP92} entre autres. 
Dans \cite{CG96}, Coates et Greenberg ont formul\'e une
conjecture assez g\'en\'erale d\'ecrivant le module des normes
universelles. La formulation du r\'esultat que nous
d\'emontrons est due \`a Nekov\'a\v{r} et concerne les normes
universelles dans l'extension cyclotomique pour une
repr\'esentation $p$-adique de de Rham. Si on l'applique au
module de Tate d'une courbe elliptique $E$, on retrouve, via la th\'eorie
de Kummer, certains r\'esultats des auteurs
cit\'es pr\'ec\'edemment: si $E$ est
ordinaire, alors  $\mathcal{N}_{K_\infty/K}(E)$ est un
$\Lambda_K$-module de rang $1$ et si $E$ est
supersinguli\`ere, alors il est nul.

\vspace{2mm}
\noindent
\textbf{Remerciements}: Je remercie Pierre Colmez pour ses nombreuses
suggestions, de la d\'emonstration du th\'eor\`eme
principal \`a la r\'edaction finale de cet article. Je remercie aussi
Jan Nekov\'a\v{r} pour ses encouragements et ses commentaires.

\renewcommand{\thesection}{\Roman{section}}

\section{Alg\`ebre diff\'erentielle des $(\phi,\Gamma)$-modules}\label{I}
L'objet de ce premier chapitre est de rappeler 
et compl\'eter certains points de la th\'eorie des 
$(\phi,\Gamma)$-modules. Ensuite, on r\'esout un probl\`eme
d'alg\`ebre diff\'erentielle. Dans le deuxi\`eme chapitre, on verra
comment cela s'applique aux repr\'esentations $p$-adiques.

\subsection{Les $(\phi,\Gamma)$-modules}\label{I1}
Dans tout cet article, $k$ d\'esigne
un corps parfait de caract{\'e}ristique $p$, 
et $K$ est une extension
finie de $F$, 
le corps des fractions de l'anneau des vecteurs de Witt sur $k$. 
On \'ecrit $\mu_{p^n} \subset \overline{K}$ pour d\'esigner 
l'ensemble des racines $p^n$-\`emes de l'unit\'e, et on d\'efinit
$K_n=K(\mu_{p^n})$ ainsi que $K_{\infty} = \cup_{n=0}^{+\infty}
K_n$. 
Soit $G_K=\on{Gal}(\overline{K}/K)$ et 
$H_K=\on{Gal}(\overline{K}/K_\infty)$ 
le noyau du caract{\`e}re cyclotomique 
$\chi: G_K \ra \Zp^*$ et $\Gamma_K=G_K/H_K$
le groupe de Galois de $K_{\infty}/K$,
qui s'identifie via le caract{\`e}re 
cyclotomique {\`a} un sous-groupe ouvert de 
$\Zp^*$. Enfin, soit $F'$ l'extension maximale 
non-ramifi\'ee de $F$
contenue dans $K_\infty$ et $k'$ le corps r\'esiduel de $F'$. 
On note $\sigma$ le Frobenius absolu (qui rel\`eve $x \mapsto
x^p$ sur $k'$).
 
D\'efinissons ici quelques anneaux de s\'eries formelles
(ces constructions sont faites en d\'etail dans \cite{PC03}): 
si $r$ est
un r\'eel positif, soit $\bdag{,r}_F$ l'anneau des s\'eries formelles
$f(X)=\sum_{k \in \ZZ} a_k X^k$ o\`u $\{a_k \in F\}_{k \in \ZZ}$ est une suite
born\'ee telle que $f(X)$ converge sur la couronne $0 < v_p(X) \leq
1/r$. Cet anneau est muni
d'une action de $\Gamma_F$, qui est triviale sur les coefficients
et donn\'ee par $\gamma(X)=(1+X)^{\chi(\gamma)}-1$ et on peut
d\'efinir un Frobenius $\phi: \bdag{,r}_F \ra \bdag{,pr}_F$ qui est
$\sigma$-semi-lin\'eaire sur les coefficients
et tel que $\phi(X)=(1+X)^p-1$. Le {\og th\'eor\`eme de 
pr\'eparation de Weierstrass \fg} montre que $\bdag{}_F = \cup_{r \geq
0} \bdag{,r}_F$ est un corps. Ce corps n'est pas complet pour la
norme de Gauss et on appelle $\mathbf{B}_F$ son
compl\'et\'e qui est un corps local de dimension $2$ dont le corps
r\'esiduel s'identifie \`a $k((\overline{X}))$.
 
L'extension $K_\infty / F_\infty$ est une extension finie de 
degr\'e de ramification $e_K \leq [K_\infty:F_\infty]$ 
et par la th\'eorie du corps de
normes de \cite{FW79,WI83} il lui correspond une extension s\'eparable
$k'((\overline{Y})) / k((\overline{X}))$ de 
degr\'e $[K_\infty:F_\infty]$ qui nous permet de d\'efinir des
extensions non-ramifi\'ees $\mathbf{B}_K / \mathbf{B}_F$ et 
$\bdag{}_K / \bdag{}_F$ de degr\'e $[K_\infty:F_\infty]$. 
On peut montrer que 
$\bdag{}_K = \cup_{r \geq 0} \bdag{,r}_K$ o\`u  
$\bdag{,r}_K$ est un $\bdag{,r}_F$-module libre de rang
$[K_\infty:F_\infty]$ qui s'identifie
\`a un anneau de s\'eries formelles
$f(Y)=\sum_{k \in \ZZ} a_k Y^k$ o\`u $\{a_k \in F'\}_{k \in \ZZ}$ est une suite
born\'ee telle que $f(Y)$ converge sur la couronne $0 < v_p(Y) \leq
1/e_K r$. L'\'el\'ement $\overline{Y}$ v\'erifie une \'equation d'Eisenstein
sur $k'((\overline{X}))$  qu'on peut relever en une \'equation sur
$\bdag{,r}_{F'}$; l'action de $\Gamma_K$
s'\'etend naturellement \`a $\bdag{,r}_K$ de m\^eme que 
le Frobenius $\phi: \bdag{,r}_K \ra \bdag{,pr}_K$.
 
Un $(\phi,\Gamma)$-module est un $\bdag{}_K$-espace vectoriel
$D^\dag$ de dimension finie, muni d'un Frobenius $\phi: D^\dag \ra 
D^\dag$ et d'une action de $\Gamma_K$ qui sont semi-lin\'eaires par
rapport \`a ceux de $\bdag{}_K$. 
On dit que $D^\dag$ est \emph{\'etale} si
$D=\mathbf{B}_K \otimes_{\bdag{}_K} 
D^\dag$ poss\`ede un r\'eseau $D_0$
stable par $\phi$ sur l'anneau des 
entiers $\mathbf{A}_K$ de $\mathbf{B}_K$, 
tel que $\phi(D_0)$ engendre $D_0$ sur
$\mathbf{A}_K$. 

D\'efinissons l'op\'erateur $\psi: D^\dag \ra D^\dag$ qui nous servira
dans la suite. On peut montrer que tout \'el\'ement $x \in D^\dag$
s'\'ecrit de mani\`ere unique $x=\sum_{i=0}^{p-1} (1+X)^i \phi(x_i)$.
\begin{defi}\label{psi}
Si $x=\sum_{i=0}^{p-1} (1+X)^i \phi(x_i)$, alors
on pose $\psi(x)=x_0$.
\end{defi}

Ceci fait de $\psi$ un inverse \`a
gauche de $\phi$ qui commute \`a l'action de $\Gamma_K$ et qui
v\'erifie $\psi(\phi(x)y)=x \psi(y)$ si $x \in \bdag{}_K$ et $y \in D^\dag$.

Il existe $r(K)$ tel que si $p^{n-1}(p-1) \geq 
r \geq r(K)$, alors on a une application injective $\iota_n :
\bdag{,r}_K \ra K_n[[t]]$ (c'est l'application $\phi^{-n}$ de
\cite[\S III.2]{CC99}). 
Par exemple si $K=F$, alors
$\iota_n(X)=\eps^{(n)} \exp(t/p^n)-1$ o\`u $\eps^{(n)}$ est une racine
primitive $p^n$-\`eme de $1$ et $\iota_n$ agit par $\sigma^{-n}$ sur les
coefficients.

On peut montrer 
(voir \cite{FC96})
que l'ensemble des 
sous $\bdag{,r}_K$-modules de type fini $M$ de
$D^\dag$ tels que $\phi(M) \subset 
\bdag{,pr}_K \otimes_{\bdag{,r}_K} M$ admet un plus
grand \'el\'ement $D^{\dag,r}$ et qu'il existe
$r(D)$ que l'on peut supposer $\geq r(K)$
tel que si $r \geq r(D)$, alors $D^\dag =
\bdag{}_K \otimes_{\bdag{,r}_K} D^{\dag,r}$.
On utilise alors l'application $\iota_n$ pour d\'efinir 
$K_n[[t]] \otimes^{\iota_n}_{\bdag{,r}_K} D^{\dag,r}$ et
$K_n((t)) \otimes^{\iota_n}_{\bdag{,r}_K} D^{\dag,r}$.

Le lemme suivant sera utile par la suite:

\begin{lemm}\label{gam_tor_pg}
Si $y \in (D^\dag)^{\psi=1}$, alors il existe $P(\gamma) \in
F'[\Gamma_K]$ tel que $P(\gamma)y=0$ si et seulement si $y \in
(D^\dag)^{\phi=1}$. 
\end{lemm}

\begin{proof}
Un petit calcul montre que $(D^\dag)^{\phi=1}$ est un $\Qp$-espace
vectoriel de dimension $\leq \dim(D^\dag)$ 
(il suffit de remarquer que des \'el\'ements de $(D^\dag)^{\phi=1}$
qui sont li\'es sur $\bdag{}_K$ le sont sur
$\Qp=(\bdag{}_K)^{\phi=1}$) et 
donc qu'il existe $P(\gamma) \in
\Qp[\Gamma_K]$ en fait
tel que $P(\gamma)$ annule $(D^\dag)^{\phi=1}$. 
Montrons donc la r\'eciproque.

Supposons que $(\sum a_i \gamma^i) y = 0$ est une relation de longueur
minimale avec $a_i \in F'$. On peut supposer que l'un des $a_i$ est
\'egal \`a $1$. En appliquant $\psi$ et en utilisant le fait
que d'une part $\psi(y)=y$ et que d'autre part $\psi$ agit par
$\sigma^{-1}$ sur $F'$, on voit que $a_i \in \Qp$ pour tout
$i$ et on suppose donc \`a partir de maintenant que $P(\gamma) \in
\Qp[\Gamma_K]$.

Nous utiliserons ci-dessous le r\'esultat suivant: si $P(\gamma) \in
\Qp[\Gamma_K]$, alors il existe une constante $C(P,d)$ telle que pour
tout $M$ qui est un $K_\infty[[t]]$-module libre de rang $d$, 
muni d'une action semi-lin\'eaire de $\Gamma_K$ par automorphismes, le
$F'$-espace vectoriel $M^{P(\gamma)=0}$ est de dimension $\leq C(P,d)$.

Fixons $r \geq r(D)$ et consid\'erons, pour $n$ tel que
$p^{n-1}(p-1) \geq r$, le $K_n[[t]]$-module libre de rang $d$ d\'efini
ci-dessus: $K_n[[t]] \otimes^{\iota_n}_{\bdag{,r}_K} D^{\dag,r}$. On voit
que l'on a une injection: 
\[ (D^{\dag,r})^{P(\gamma)=0} \hookrightarrow \left(
K_\infty[[t]] \otimes_{K_n[[t]]} K_n[[t]]
\otimes^{\iota_n}_{\bdag{,r}_K} D^{\dag,r} \right)^{P(\gamma)=0}, \]
ce qui fait que $(D^{\dag,r})^{P(\gamma)=0}$ est un
$F'$-espace vectoriel de dimension $\leq C(P,d)$ et donc que
$(D^{\dag})^{P(\gamma)=0} = \cup_{r \geq r(D)} (D^{\dag,r})^{P(\gamma)=0}$ est un
$F'$-espace vectoriel de dimension $\leq C(P,d)$. Comme $\phi$ commute
\`a $P(\gamma)$, $(D^{\dag})^{P(\gamma)=0}$ est un
$F'$-espace vectoriel de dimension finie et stable par $\phi$ ce qui
fait que $\phi:(D^{\dag})^{P(\gamma)=0} \ra (D^{\dag})^{P(\gamma)=0}$ 
est bijectif.

Si $z \in (D^{\dag})^{\psi=0,P(\gamma)=0}$, alors $z
= \phi(w)$ pour un $w \in (D^\dag)^{P(\gamma)=0}$ et donc 
$0=\psi(z)=w$ ce qui fait que $z=0$ et donc que
$(D^\dag)^{\psi=0,P(\gamma)=0}=0$
(ceci g\'en\'eralise un r\'esultat de \cite{CC98}). 
Pour conclure, il suffit de
remarquer que si $\psi(y)=y$ et $P(\gamma)(y)=0$, 
alors $\psi(1-\phi)y=0$
et donc $y=\phi(y)$. 
\end{proof}

\subsection{$(\phi,\Gamma)$-modules de de Rham}\label{secpgdr}
Dans ce paragraphe, on d\'efinit les $(\phi,\Gamma)$-modules de de
Rham et on rappelle certains des r\'esultats de \cite[\S 5]{LB02} \`a
leur sujet.

\begin{defi}\label{pg_dr}
On dit qu'un $(\phi,\Gamma)$-module $D^\dag$ est de \emph{de Rham}, si
et seulement s'il existe $r\in\RR$, $n\in\NN$ avec $p^{n-1}(p-1) \geq r \geq
r(D)$, tels que le $K$-espace vectoriel 
\[ \left(K_n((t)) \otimes^{\iota_n}_{\bdag{,r}_K} D^{\dag,r}\right)^{\Gamma_K} \] 
est de dimension $d=\dim(D^\dag)$.
\end{defi}

Si c'est le cas, alors $(K_n((t)) \otimes^{\iota_n}_{\bdag{,r}_K} D^{\dag,r})^{\Gamma_K}$ 
est de dimension $d=\dim(D^\dag)$ pour tous les 
$r\in\RR$, $n\in\NN$ tels que $p^{n-1}(p-1) \geq r \geq r(D)$.

Nous allons rappeler les r\'esultats de \cite[\S 4,5]{LB02} qui
permettent de donner une autre caract\'erisation des
$(\phi,\Gamma)$-modules de de Rham. Ces r\'esultats sont aussi 
expliqu\'es dans le {\og s\'eminaire Bourbaki \fg} \cite{PC01}.
 
L'anneau $\bdag{,r}_K$ s'identifiant \`a un anneau de s\'eries
formelles convergeant sur une couronne, 
il est naturellement muni d'une topologie de Fr\'echet, 
la topologie de la {\og convergence compacte \fg}, et
son compl\'et\'e $\bnrig{,r}{,K}$ pour cette topologie s'identifie \`a 
l'anneau de s\'eries formelles
$f(Y)=\sum_{k \in \ZZ} a_k Y^k$ o\`u $\{a_k \in 
F'\}_{k \in \ZZ}$ est une suite
non n\'ecessairement born\'ee telle que 
$f(Y)$ converge sur la couronne $0 < v_p(Y) \leq 1/e_K r$.
Par exemple, si on pose $t=\log(1+X)$, alors $t \in \bnrig{,r}{,F} \subset 
\bnrig{,r}{,K}$ pour tout $r \geq 0$.
L'anneau $\bnrig{}{,K} = \cup_{r \geq 0} \bnrig{,r}{,K}$ est 
{\og l'anneau de Robba \fg}. Ces anneaux ont \'et\'e \'etudi\'es dans
\cite[\S 4]{LB02} et nous allons rappeler quelques-uns des r\'esultats
qui nous seront utiles dans la suite.

L'application $\iota_n$
se prolonge en une application injective 
$\iota_n : \bnrig{,r}{,K} \ra K_n[[t]]$. 
L'action de $\Gamma_K$ sur $\bnrig{,r}{,K}$ s'\'etend en une action de
l'alg\`ebre de Lie de $\Gamma_K$ donn\'ee par $\nabla(f) =
\log(\gamma)(f)/\log_p(\chi(\gamma))$ pour $\gamma \in \Gamma_K$ assez
proche de $1$. 
Si $f = f(X) \in \bnrig{,r}{,F}$ alors
$\nabla(f(X)) = t(1+X)df/dX$. Si $f \in \bnrig{,r}{,K}$ alors
on pose $\partial(f) = t^{-1} \nabla(f)$ ce qui fait que si 
$f = f(X) \in \bnrig{,r}{,F}$ alors
$\partial(f(X)) = (1+X)df/dX$ et que si $f \in \bnrig{}{,K}$ v\'erifie
une \'equation alg\'ebrique $P(f)=0$ sur $\bnrig{}{,F}$ 
telle que $P'(f) \neq 0$, alors
on peut aussi calculer $\partial(f)$ par la formule
$\partial(f)= -(\partial P)(f)/P'(f)$. En particulier 
$\partial(f)=0$ si et seulement si $f \in F'$.
 
Si $D^\dag$ est un $(\phi,\Gamma)$-module, on d\'efinit $\ddagrig$ par
$\ddagrig = \bnrig{}{,K} \otimes_{\bdag{}_K} D^\dag$. 
L'alg\`ebre de Lie de $\Gamma_K$ agit sur $\ddagrig$ par la
formule $\nabla_D(x) =
\log(\gamma)(x)/\log_p(\chi(\gamma))$ pour $\gamma \in \Gamma_K$ assez
proche de $1$ et on a donc aussi une application $\partial_D =
t^{-1} \nabla_D : \ddagrig \ra t^{-1} \ddagrig$.  La proposition
suivante se d\'emontre exactement de la m\^eme mani\`ere que 
\cite[th\'eor\`eme 5.10]{LB02}.

\begin{prop}\label{ndr_pg}
Si $D^\dag$ est de de Rham, si $r \geq r(D)$, et si
$N_r$ est l'ensemble des $x \in D^{\dag,r}_{\mathrm{rig}}[1/t]$
tels que pour tout
$n$ tel que $p^{n-1}(p-1) \geq r$, on ait
\[ \iota_n(x) \in K_n[[t]] \otimes_K   
(K_n((t)) \otimes^{\iota_n}_{\bdag{,r}_K} D^{\dag,r})^{\Gamma_K}, \] 
et si on pose $N = \bnrig{}{,K} \otimes_{\bnrig{,r}{,K}} N_r$,
alors $N$ est un $\bnrig{}{,K}$-module libre de rang $d$ stable par les
actions induites de $\phi$ et de $\Gamma_K$, tel que $\partial_D(N)
\subset N$ et tel que $N[1/t]=\ddagrig[1/t]$.
\end{prop}

Nous \'etablirons quelques propri\'et\'es de $N$ dans le paragraphe \ref{II1}.

\begin{rema}
Donnons nous une variable $\ell_Y$ sur laquelle on fait agir $\gamma \in
\Gamma_K$ par $\gamma(\ell_Y)=\ell_Y + \log(\gamma(Y)/Y)$.
Si $d=\dim(D^\dag)$, alors on dit que le 
$(\phi,\Gamma)$-module $D^\dag$ est
\begin{enumerate}
\item \emph{cristallin}, si $(\ddagrig[1/t])^{\Gamma_K}$ est un $F$-espace
  vectoriel de dimension $d$;
\item \emph{semi-stable}, si $(\ddagrig[\ell_Y][1/t])^{\Gamma_K}$ est un $F$-espace
  vectoriel de dimension $d$.
\end{enumerate}

Il n'est pas difficile de voir que {\og cristallin \fg} 
implique {\og semi-stable \fg}
et que {\og semi-stable \fg} 
implique {\og de de Rham \fg}; on a d'ailleurs dans ce cas
\[ N=\left( \bnrig{}{,K}[\ell_Y] \otimes_F
(\ddagrig[\ell_Y][1/t])^{\Gamma_K} \right)^{d/d\ell_Y=0}. \] 
De plus, le th\'eor\`eme de monodromie $p$-adique pour les \'equations
diff\'erentielles (d\'emontr\'e ind\'ependamment 
dans \cite{YA01,KK00,ZM01}) montre que $D^\dag$ est de de Rham
si et seulement s'il existe une extension finie $L$ de $K$ telle que $\bdag{}_L
\otimes_{\bdag{}_K} D^\dag$ est semi-stable. Nous n'utiliserons pas
cela par la suite.
\end{rema}

\subsection{Alg\`ebre diff\'erentielle}\label{I2}
On aura besoin dans la suite d'un argument qui est une variante du
{\og d\'eterminant Wronskien \fg} et qui est l'objet de ce chapitre.
Soient $H$ un corps diff\'erentiel, dont on notera $\partial$ la
d\'erivation, et $k$, $s$ et $v$ trois entiers $\geq 1$. 
On \'etend naturellement $\partial$ \`a $\partial: H^v \ra H^v$. 
On se donne $s+1$ vecteurs $x_1,\cdots,x_{s+1} \in H^v$ qui 
satisfont les deux conditions ci-dessous:
\begin{enumerate}
\item d'une part, les
$s$ vecteurs $X^{k-1}_w=(x_w,\partial x_w, \cdots, \partial^{k-1} x_w)$
pour $1 \leq w \leq s$ sont lin\'eairement ind\'ependants sur $H$;
\item d'autre part, les $s+1$ vecteurs 
$X^k_w=(x_w,\partial x_w, \cdots, \partial^k x_w)$
pour $1 \leq w \leq s+1$ sont lin\'eairement d\'ependants sur $H$.
\end{enumerate}

\begin{prop}\label{wronski}
Sous les hypoth\`eses ci-dessus, 
les $s+1$ vecteurs $\{X^k_w\}$ 
pour $1 \leq w \leq s+1$ sont li\'es
sur $H^{\partial=0}$. 
\end{prop}

\begin{proof}
Comme on suppose que les $X^k_w$ sont li\'es sur $H$, il existe des
\'el\'ements $\lambda_1,\cdots,\lambda_{s+1} \in H$ tels que
$\sum_{w=1}^{s+1} \lambda_w X^k_w=0$, ce qui se traduit par le fait
que pour tout $0 \leq i \leq k$, on
a la relation $R_i$:
\[\sum_{w=1}^{s+1} \lambda_w \partial^i (x_w) = 0. \]
Comme on a suppos\'e que les $\{X^{k-1}_w\}_{1 \leq w \leq s}$
sont libres, on voit que $\lambda_{s+1} \neq 0$ et on peut donc
supposer que $\lambda_{s+1} = 1$ ce que l'on fait \`a partir de
maintenant.

Si on d\'erive $R_i$, on trouve que
\[\sum_{w=1}^{s+1} \lambda_w \partial^{i+1} (x_w) + 
\sum_{w=1}^{s+1} \partial
(\lambda_w) \partial^i (x_w) = 0. \]
Si $0 \leq i \leq k-1$, alors le premier terme de la relation
ci-dessus est nul par $R_{i+1}$ et on trouve que
\[ \sum_{w=1}^{s+1} \partial (\lambda_w) \partial^i (x_w) = 0 \]
pour $0 \leq i \leq k-1$ et donc que
\[ \sum_{w=1}^{s} \partial (\lambda_w) \partial^i (x_w) = 0 \]
puisque $\lambda_{s+1} = 1$. Ceci nous donne une relation 
\[ \sum_{w=1}^{s} \partial (\lambda_w) X^{k-1}_w = 0 \]
entre les 
$\{X^{k-1}_w\}_{1 \leq w \leq s}$ et comme on a suppos\'e que
ceux-ci sont libres, c'est que $\partial(\lambda_w)=0$ pour tout 
$1 \leq w \leq s+1$. 

On a donc montr\'e que si les $\{X^k_w\}_{1 \leq w \leq s+1}$ sont li\'es
sur $H$, alors ils sont li\'es sur $H^{\partial=0}$. 
\end{proof}

En consid\'erant la premi\`ere {\og composante \fg}
de $X^k_w$, on trouve:

\begin{coro}\label{wronski_coro}
Sous les hypoth\`eses ci-dessus,
il existe des constantes $\lambda_w
\in H^{\partial=0}$ pour $1 \leq w \leq s+1$, non toutes nulles,  
telles que \[ \sum_{w=1}^{s+1} \lambda_w x_w = 0. \]
\end{coro}

\subsection{Normes universelles: $(\phi,\Gamma)$-modules positifs}\label{I3}
Dans ce paragraphe, nous allons montrer la proposition suivante:

\begin{prop}\label{pg_diff}
Soit $D^\dag$ un
$(\phi,\Gamma)$-module de de Rham qui v\'erifie $D^\dag \subset N$
et $D^\dag \cap t N = 0$. Si $y \in D^\dag$ et $k \geq 0$ sont tels
que $\partial_D^k y \in t N$, alors il existe $P(\gamma) \in
F'[\Gamma_K]$ tel que $P(\gamma)y=0$. 
\end{prop}

\begin{proof}
Si $k=0$, alors le fait que $D^\dag \cap t N = 0$ implique
imm\'ediatement que $y=0$. On suppose d\'esormais 
que $k \geq 1$, et on se fixe $\gamma \in
\Gamma_K$ d'ordre infini. Comme $D^\dag$ est un $\bdag{}_K$-espace
vectoriel de dimension finie, 
il existe un entier $v \geq 1$ tel
que $y, \gamma(y), \cdots, \gamma^{v-1}(y)$ sont libres 
sur  $\bdag{}_K$, mais 
$y, \gamma(y), \cdots, \gamma^v(y)$ sont li\'es sur  $\bdag{}_K$, ce
qui fait que pour tout $w \geq v$, il existe des \'el\'ements
$a_0^w,\cdots,a_{v-1}^w$ de  $\bdag{}_K$ tels que
\[ \gamma^w(y) = 
\sum_{j=0}^{v-1} a_j^w \gamma^j(y) =
a_0^w y + \cdots + a_{v-1}^w \gamma^{v-1}(y). \]

Si l'on d\'erive $k$ fois la relation ci-dessus, on trouve que:
\[ \partial_D^k( \gamma^w(y) ) =
\sum_{j=0}^{v-1} \sum_{i=0}^{k-1} \binom{k}{i} \partial^i(a_j^w)
\partial_D^{k-i}(\gamma^j(y)) + \sum_{j=0}^{v-1} \partial^k(a_j^w) 
\gamma^j(y). \]  
Supposons que l'on se donne un entier $s \geq 1$ et des \'el\'ements
$\mu_v,\cdots,\mu_{v+s}$ de $\bdag{}_K$ tels que pour tout $0 \leq i
\leq k-1$ et pour tout $0 \leq j \leq v-1$, on ait $\sum_{w=v}^{v+s}
\mu_w \partial^i(a_j^w)=0$. Bien s\^ur, de tels \'el\'ements (non tous
nuls) existent si $s \gg 0$. 
On voit alors que
\begin{equation}\label{eqn01}
\sum_{w=v}^{v+s} \mu_w  \partial_D^k( \gamma^w(y) ) = 
  \sum_{j=0}^{v-1} \left( \sum_{w=v}^{v+s} \mu_w \partial^k(a_j^w) \right)
\gamma^j(y). 
\end{equation}
Montrons que cela implique que $\sum_{w=v}^{v+s}
\mu_w \partial^i(a_j^w)=0$ pour $i=k$ et pour $0 \leq j \leq v-1$, ce
qui est \'equivalent au fait que les deux termes de l'\'equation (\ref{eqn01})
ci-dessus sont nuls.

On a $\partial_D^k( \gamma^w(y) ) = \chi(\gamma)^{kw}
\gamma^w (\partial_D^k (y)) \in t N$ et le terme de gauche de
l'\'equation (\ref{eqn01}) est donc dans $t N$ tandis que le terme de
droite est dans $D^\dag$ ce qui fait que les deux termes sont nuls
puisqu'on a suppos\'e que $D^\dag \cap t N = 0$. 

Pour $w \geq 1$, posons $x_w=(a_0^{w+v-1},\cdots,a_{v-1}^{w+v-1})$ et
$X^{k-1}_w=(x_w,\cdots,\partial^{k-1} x_w)$ ainsi que 
$X^k_w=(x_w,\cdots,\partial^k x_w)$.
Si $s \geq 1$ est le plus petit entier tel
que les $X^{k-1}_w$ pour $1 \leq w \leq s$ sont libres et les 
$X^{k-1}_w$ pour $1 \leq w \leq s+1$ sont li\'es, alors les calculs
pr\'ec\'edents montrent \emph{qu'en plus}, les 
$X^k_w$ pour $1 \leq w \leq s+1$ sont \emph{eux aussi} li\'es. 
On peut alors appliquer le corollaire
\ref{wronski_coro} pour en d\'eduire l'existence de $\lambda_w \in
(\bdag{}_K)^{\partial=0}=F'$ pour $1 \leq w \leq s+1$ tels que 
\[ \sum_{w=1}^{s+1} \lambda_w x_w = 0. \]
Si l'on combine cela avec le fait que par d\'efinition, on a 
$\gamma^w(y) = 
a_0^w y + \cdots + a_{v-1}^w \gamma^{v-1}(y)$,
on trouve que
\[ \sum_{w=1}^{s+1} \lambda_w \gamma^{w+v-1}(y) = 0. \]
Ceci montre bien qu'il existe 
$P(\gamma) \in F'[\Gamma_K]$ tel que $P(\gamma)(y)=0$.
\end{proof}

\section{Repr\'esentations $p$-adiques et normes universelles}\label{II}
L'objet de ce chapitre est de montrer comment le r\'esultat
d\'emontr\'e au chapitre pr\'ec\'edent (la proposition \ref{pg_diff})
implique le th\'eor\`eme principal de cet article. Pour cela, on
rappelle la correspondance entre repr\'esentations $p$-adiques et
$(\phi,\Gamma)$-modules, puis on utilise la formule r\'eciprocit\'e de
Cherbonnier-Colmez pour se ramener \`a la situation du chapitre
pr\'ec\'edent. On conclut par un argument de {\og d\'evissage \fg}.

\subsection{Repr\'esentations $p$-adiques et $(\phi,\Gamma)$-modules}\label{II1}
Une repr\'esentation $p$-adique est un $\Qp$-espace vectoriel $V$
de dimension finie $d=\dim_{\Qp}(V)$, muni d'une action lin\'eaire et
continue de $G_K$. Afin d'\'etudier les repr\'esentations
$p$-adiques, Fontaine a construit 
(voir \cite{Bu88per} par exemple) un
certains nombre d'anneaux de p\'eriodes $p$-adiques, ce qui conduit
\`a la d\'efinition des 
repr\'esentations \emph{cristallines}, \emph{semi-stables}
ou de \emph{de Rham}.
 
D'autre part, en combinant les constructions de Fontaine (voir
\cite{FO91}) et un th\'eor\`eme de Cherbonnier-Colmez (voir
\cite{CC98}), 
on d\'efinit un foncteur $V \mapsto \ddag{}(V)$  
qui induit une \'equivalence de cat\'egories entre la cat\'egorie
des repr\'esentations $p$-adiques et la cat\'egorie des
$(\phi,\Gamma)$-modules \'etales. 
%On note $\mathbf{V}$ l'inverse de ce foncteur. 
Si par exemple $H_K$ agit trivialement sur $V$, alors $\ddag{}(V)=
\bdag{}_K \otimes_{\Qp} V$ et on r\'ecup\`ere $V$ \`a partir de
$\ddag{}(V)$ par $V = \ddag{}(V)^{\phi=1}$. En g\'en\'eral, 
la situation est plus compliqu\'ee mais on peut
quand m\^eme montrer que $\ddag{}(V)^{\phi=1}=V^{H_K}$ ce dont nous
aurons besoin par la suite.

Rappelons que pour un $(\phi,\Gamma)$-module $D^\dag$, 
on a d\'efini au paragraphe \ref{I1} un r\'eel $r(D)$ tel que si $r
\geq r(D)$, alors $D^\dag = \bdag{}_K \otimes_{\bdag{,r}_K}
D^{\dag,r}$ et on posera dans la suite $r(V)=r(\ddag{}(V))$.
Les r\'esultats de \cite{CC99} montrent que si $V$ est une
repr\'esentation $p$-adique, et si $p^{n-1}(p-1) \geq r \geq r(V)$, 
alors $\ddR(V)=(K_n((t)) \otimes^{\iota_n}_{\bdag{,r}_K} \ddag{,r}(V))^{\Gamma_K}$ et
donc que $V$ est de de Rham si et seulement si
$\ddag{}(V)$ est de de Rham au sens de la d\'efinition \ref{pg_dr}. 
On dispose alors d'une part des applications 
\[ \iota_n : \ddag{,r}(V) \ra K_n((t)) \otimes_K \ddR(V) \] et d'autre
part de l'\'equation diff\'erentielle $p$-adique $\ndr(V)$ qui est le
module $N$ dont on a rappel\'e la construction dans la proposition
\ref{ndr_pg}. 

\begin{rema}
Les r\'esultats de
\cite{LB02} montrent par ailleurs
qu'une repr\'esentation $V$ est cristalline (ou
semi-stable) si et seulement si son 
$(\phi,\Gamma)$-module $\ddag{}(V)$ est cristallin (ou
semi-stable). On retrouve alors les invariants
associ\'es \`a $V$ par la th\'eorie de Hodge $p$-adique de la
mani\`ere suivante: $\dcris(V) = (\dnrig{}(V)[1/t])^{\Gamma_K}$ et 
$\dst(V) = (\dnrig{}(V)[\ell_Y][1/t])^{\Gamma_K}$. 
\end{rema}
 
Nous allons avoir besoin de quelques r\'esultats concernant
$\ndr(V)$, et on suppose \`a partir de maintenant que $V$ est de de Rham.

\begin{lemm}\label{HT_ndr}
Si $V$ est une repr\'esentation de de Rham, 
dont les poids de Hodge-Tate sont $\geq 0$, alors $\ddag{}(V) \subset
\ndr(V)$ et si  $\ddag{}(V) \subset t \ndr(V)$, alors 
les poids de Hodge-Tate de $V$ sont $\geq 1$.
\end{lemm}

\begin{proof}
\'Etant donn\'ee la construction de $N=\ndr(V)$ que l'on a donn\'ee
dans la proposition \ref{ndr_pg}, on voit qu'il suffit 
pour montrer le premier point de montrer  
que si $r \geq r(V)$ et 
$y \in \ddag{,r}(V)$ et si $p^{n-1}(p-1) \geq r$, 
alors $\iota_n(y) \in K_n[[t]] \otimes_K \ddR(V)$. On sait 
par les constructions de \cite[\S 5]{LB02} que $\iota_n(y) \in \bdR^+
\otimes_{\Qp} V$ et le r\'esultat suit du fait que si les poids de
Hodge-Tate de $V$ sont $\geq 0$, alors 
$\bdR^+ \otimes_{\Qp} V \subset \bdR^+ \otimes_K \ddR(V)$. 

D'autre part, la d\'emonstration de 
\cite[proposition 5.15]{LB02} montre que
le $\bdR^+$-module engendr\'e par $\iota_n(N_r)$ est $\bdR^+ \otimes_K
\ddR(V)$. De m\^eme, le $\bdR^+$-module engendr\'e par
$\iota_n(\ddag{,r}(V))$ est $\bdR^+ \otimes_{\Qp} V$ et
la deuxi\`eme assertion suit du fait que si 
$\bdR^+ \otimes_{\Qp} V \subset t \bdR^+ \otimes_K \ddR(V)$, alors les
poids de Hodge-Tate de $V$ sont $\geq 1$.
\end{proof}

\begin{lemm}\label{ndr_sub}
Si $V$ est une repr\'esentation de de Rham dont les poids de
Hodge-Tate sont $\geq
0$ et si $W$ est une sous-repr\'esentation de $V$, alors $\ndr(W) =
\dnrig{}(W) \cap \ndr(V)$.
\end{lemm}

\begin{proof}
On voit que $\dnrig{}(W) \cap \ndr(V)$ est un sous
$\bnrig{}{,K}$-module de $\dnrig{}(W)$ qui est libre de rang $\dim(W)$ 
(parce que $\ndr(V)[1/t]=\dnrig{}(V)[1/t]$) et
qui est stable par $t^{-1} \nabla_W$ (qui co\"{\i}ncide avec 
l'op\'erateur induit par $t^{-1} \nabla_V$),
ce qui fait que $\ndr(W) = \dnrig{}(W) \cap \ndr(V)$ 
par d\'efinition de $\ndr(W)$ (voir \cite[th\'eor\`eme 5.10]{LB02}).
\end{proof}

\begin{lemm}\label{no_divis}
Si $V$ est une repr\'esentation de de Rham dont les poids de
Hodge-Tate sont $\geq 0$ et qui  
n'admet pas de sous-repr\'esentation dont les poids de
Hodge-Tate sont $\geq 1$, alors 
$\ddag{}(V) \cap t \ndr(V) = \{0\}$.
\end{lemm}

\begin{proof}
Si $z \in \ddag{}(V) \cap t\ndr(V)$, alors on voit que l'intersection
$\ddag{}(z)$ des sous $\bdag{}_K$-espaces vectoriels de $\ddag{}(V)$
stables par $\phi$ et $\Gamma_K$ 
et contenant $z$
v\'erifie $\ddag{}(z) \subset
t\ndr(V)$. La proposition \cite[1.1.6]{FO91} implique 
d'autre part que $\ddag{}(z)$
est un sous $(\phi,\Gamma)$-module \'etale de $\ddag{}(V)$ et donc
qu'il existe une repr\'esentation $p$-adique $V_z \subset V$ telle que
$\ddag{}(z) = \ddag{}(V_z)$. On voit alors que 
$\ddag{}(V_z) \subset  t\ndr(V) \cap \dnrig{}(V_z) = t\ndr(V_z)$ 
par le lemme \ref{ndr_sub} et
donc que les poids de Hodge-Tate de $V_z$ sont tous $\geq 1$, 
par la proposition \ref{HT_ndr}, et donc
que $V_z=0$ ce qui fait que $z=0$.
\end{proof}

Rappelons que l'on a une application
$\iota_n : \ddag{,r}(V) \ra K_n((t)) \otimes_K \ddR(V)$ et on en d\'eduit une
application \[ \delta_{V(-k)} \circ \iota_n : \ddag{,r}(V) \ra
K_n \otimes_K \ddR(V) \] o\`u $\delta_{V(-k)}: K_n((t)) 
\otimes_K \ddR(V) \ra K_n \otimes_K \ddR(V)$ est l'application 
{\og coefficient de $t^k$ \fg}.

\begin{lemm}\label{divis_t}
Si $z \in \ndr(V)$ est tel que pour tout $n \gg 0$, on ait $\iota_n(z)
\in t K_n[[t]] \otimes_K \ddR(V)$, alors $z \in t \ndr(V)$. En
particulier, si $\delta_V \circ \iota_n(z)=0$  pour tout $n \gg 0$,
alors $z \in t \ndr(V)$.
\end{lemm}

\begin{proof}
Fixons $r \gg 0$ tel que $z \in N_r$.
La d\'emonstration de la proposition \cite[5.15]{LB02} montre
\footnote{il est en fait incorrectement affirm\'e que $\iota_n$ est
surjective - elle est surjective modulo $t^w$ pour tout $w
\geq 0$.} qu'il existe une base
$f_1, \cdots, f_d$ de $N_r$ dont les images par $\iota_n$ sont une
base de $K_n[[t]] \otimes_K \ddR(V)$ pour tout $n \gg 0$. Si l'on
\'ecrit que $z=\sum_{i=1}^d z_i f_i$, et si $\iota_n(z)
\in t K_n[[t]] \otimes_K \ddR(V)$, alors $z_i(\eps^{(n)}-1)=0$ pour
tout $n \gg 0$ et donc $t$ divise $z_i$ dans $\bnrig{}{,K}$ ce qui
fait que  $z \in t \ndr(V)$.
\end{proof}

\subsection{Cohomologie galoisienne et repr\'esentations de de Rham}\label{II2}
Rappelons que Herr a montr\'e dans \cite{LH98} comment construire les
groupes de cohomologie galoisienne $H^i(K,V)$ \`a partir du 
$(\phi,\Gamma)$-module associ\'e \`a $V$. Nous utiliserons la version
qu'en ont donn\'ee Cherbonnier et Colmez dans \cite{CC99} et que nous
rappelons bri\`evement. 

Rappelons que l'on a d\'efini (voir d\'efinition \ref{psi}) un
op\'erateur $\psi: \ddag{}(V) \ra \ddag{}(V)$.
Dans \cite[\S I.5]{CC99}, une application $h^1_{K_n,V} :
\ddag{}(V)^{\psi=1} \ra H^1(K_n,V)$ est construite pour tout $n \geq
0$. Ces applications $h^1_{K_n,V}$ donnent lieu 
(par un th\'eor\`eme non publi\'e de Fontaine dont on trouvera une
d\'emonstration dans \cite[\S II.1]{CC99})
\`a un isomorphisme entre
$\ddag{}(V)^{\psi=1}$ et la cohomologie d'Iwasawa $H^1_{\Iw}(K,V)$ par
passage \`a la limite sur $n$. 

On suppose \`a pr\'esent que $V$ est de de Rham, et on va voir \`a
quelle condition $h^1_{K_n,V}(y) \in  H^1_g(K_n,V)$ si  
$y \in \ddag{}(V)^{\psi=1}$. L'op\'erateur $\partial_D$ que l'on a
d\'efini au paragraphe \ref{secpgdr} v\'erifie $\iota_n \circ \partial_D  =
p^{-n} \partial_V \circ \iota_n$ o\`u 
\[ \partial_V = d/dt \otimes \on{Id} :
 K_n((t))  \otimes_K \ddR(V) \ra K_n((t)) \otimes_K \ddR(V). \]

Si l'on combine le th\'eor\`eme 
\cite[IV.2.1]{CC99} avec des propri\'et\'es bien connues de 
{\og l'exponentielle duale de Bloch-Kato \fg}, et le fait que 
$\delta_{V(-k)} = k!^{-1} \delta_V \circ 
\partial_V^k$ pour $k \geq 0$, alors on trouve que:
\begin{prop}\label{ch_co}
Il existe $r_{\psi}(V)$ que l'on peut supposer $\geq r(V)$
tel que si $r \geq r_{\psi}(V)$,
alors $\ddag{}(V)^{\psi=1} \subset \ddag{,r}(V)$ et  
si $y \in \ddag{}(V)^{\psi=1}$, 
$p^{n-1}(p-1) \geq r \geq r_{\psi}(V)$ et $k \geq 0$, alors  
on a $h^1_{K_n,V(-k)}(y(-k)) \in  H^1_g(K_n,V(-k))$ si et seulement si 
$\delta_V \circ \partial_V^k \circ \iota_n (y) = 0$.
\end{prop}

Enfin pour terminer, on calcule la $F'[\gamma_K]$-torsion de
$\ddag{}(V)^{\psi=1}$.
\begin{lemm}\label{gam_tor}
Si $y \in \ddag{}(V)^{\psi=1}$, alors il existe $P(\gamma) \in
F'[\Gamma_K]$ tel que $P(\gamma)y=0$ si et seulement si $y \in
V^{H_K}$. 
\end{lemm}

\begin{proof}
\'Etant donn\'e que $\ddag{}(V)^{\phi=1}=V^{H_K}$ comme on l'a
rappel\'e au paragraphe \ref{II1}, c'est une cons\'equence imm\'ediate
du lemme \ref{gam_tor_pg}.
\end{proof}

\begin{rema}
Perrin-Riou a en fait d\'etermin\'e la structure du $\Qp \otimes_{\Zp}
\Lambda_K$-module $H^1_{\Iw}(K,V)$. Son sous-module de torsion
s'identifie \`a $V^{H_K}$ et $H^1_{\Iw}(K,V)/V^{H_K}$ est (au moins si
$K/\Qp$ est finie) un $\Qp \otimes_{\Zp}
\Lambda_K$-module libre de rang $[K:\Qp]\dim(V)$.
\end{rema}

\subsection{Normes universelles}\label{II3}
Dans ce paragraphe, on montre le th\'eor\`eme \ref{main_theo}. On
commence par un r\'esultat un peu plus faible (la proposition
\ref{univ_irred} ci-dessous) puis on montre comment en d\'eduire le cas
g\'en\'eral en d\'evissant un peu.

Si $W$ est une repr\'esentation $p$-adique, soit
$\ddag{}(W)^{\psi=1}_g$ l'ensemble des  
$z \in \ddag{}(W)^{\psi=1}$ tels que pour tout $n \geq 0$, 
on ait $h^1_{K_n,W}(z) \in H^1_g(K_n,W)$. 

\begin{prop}\label{univ_irred}
Si $k \geq 0$ et si $V$ est une 
repr\'esentation de de Rham dont les poids de
Hodge-Tate sont $\geq -k$ et qui  
n'admet pas de sous-repr\'esentation dont les poids de
Hodge-Tate sont $\geq 1-k$, 
alors $\ddag{}(V)^{\psi=1}_g \subset V^{H_K}$.

D'autre part, on a $\ddag{}(V)^{\psi=1}_g = 0$ si $k = 0$
et $\ddag{}(V)^{\psi=1}_g = V^{H_K}$ si $k \geq 1$.
\end{prop}

\begin{proof}
Posons $D^\dag = \ddag{}(V(k))$ et $N=\ndr(V(k))$. Comme les poids de
Hodge-Tate de $V(k)$ sont $\geq 0$, le lemme \ref{HT_ndr} montre que
$\ddagrig \subset N$. D'autre part, le fait que $V$ 
n'admet pas de sous-repr\'esentation dont les poids de
Hodge-Tate sont $\geq 1-k$ et le lemme \ref{no_divis} montrent que
$D^\dag \cap t N=0$. Enfin, la proposition \ref{ch_co} montre que si
$h^1_{K_n,V}(y) \in H^1_g(K_n,V)$ pour tout $n \geq 0$, alors 
$\delta_{V(k)} \circ \partial_{V(k)}^k \circ \iota_n (y) = 0$ 
pour tout $n \gg 0$ et 
donc que $\delta_{V(k)} \circ \iota_n \circ \partial_D^k (y) = 0$
pour tout $n \gg 0$
et le lemme \ref{divis_t} 
(appliqu\'e \`a $V(k)$)
montre que $\partial_D^k (y) \in t N$. 

On est donc en mesure d'appliquer la proposition \ref{pg_diff}, qui
montre qu'il existe $P(\gamma) \in F'[\Gamma_K]$ tel que $P(\gamma)y=0$. 
Le lemme \ref{gam_tor} montre enfin que $y \in V^{H_K}$, ce qui
d\'emontre le premier point.

Ensuite, si $k=0$ et $y \in \ddag{}(V)^{\psi=1}_g$,
alors $y \in t N \cap D^\dag = 0$ et donc
$\ddag{}(V)^{\psi=1}_g = 0$. Enfin si $k \geq 1$ et $y \in V^{H_K}$,
alors $h^1_{K_n,V}(y)=0$ pour tout $n \geq 0$ et donc 
$\ddag{}(V)^{\psi=1}_g = V^{H_K}$.
\end{proof}

Pour terminer, on donne les arguments de d\'evissage permettant de
d\'eduire le th\'eor\`eme principal 
de la proposition \ref{univ_irred}.

\begin{lemm}\label{g_ext}
Si $0 \ra V' \ra V \ra V'' \ra 0$ est une suite exacte de
repr\'esentations $p$-adiques, alors on a une suite exacte:
\[ 0 \ra \ddag{}(V')^{\psi=1}_g \ra \ddag{}(V)^{\psi=1}_g \ra
\ddag{}(V'')^{\psi=1}_g \]
\end{lemm}

\begin{proof}
La th\'eorie des $(\phi,\Gamma)$-modules nous donne une suite exacte:
\[ 0 \ra \ddag{}(V') \ra \ddag{}(V) \ra
\ddag{}(V'') \ra 0, \]
et le lemme du serpent implique que l'on a:
\[ 0 \ra \ddag{}(V')^{\psi=1} \ra \ddag{}(V)^{\psi=1} \ra
\ddag{}(V'')^{\psi=1}. \]
Le lemme r\'esulte alors du fait 
(voir \cite[\S 3.8, (iii)]{Bu88sst})
que les sous-quotients des
repr\'esentations de Rham sont de de Rham.
\end{proof}

Si $W$ est une repr\'esentation de Hodge-Tate et $j \in \ZZ$, soit
$\Fil^j W$ la plus grande sous-repr\'esentation de $W$ dont tous les
poids sont $\geq j$. 

\begin{lemm}\label{no_fil_fil}
On a $\Fil^j(W/\Fil^j W)=0$. 
\end{lemm}

\begin{proof}
Soit $f:W \ra W/\Fil^j W$ la projection naturelle. Si $X
\subset W / \Fil^j W$ a tous ses poids $\geq j$, alors on
a une suite exacte
\[ 0 \ra  \Fil^j W \ra f^{-1}(X) \ra X \ra 0 \] et donc 
$f^{-1}(X)$ est une repr\'esentation de Hodge-Tate, extension de
deux repr\'esentations de Hodge-Tate \`a poids $\geq j$, et donc
elle-m\^eme \`a poids $\geq j$, ce qui fait que
$f^{-1}(X)=\Fil^j W$ et donc que $X=0$.
\end{proof}
 
Nous pouvons enfin montrer le r\'esultat principal de cet
article. Pour cela, rappelons que l'on a un isomorphisme
$H^1_{\Iw}(K,V) = \ddag{}(V)^{\psi=1}$.

\begin{theo}\label{main_in}
Si $V$ est une repr\'esentation de de Rham, alors
\begin{enumerate}
%\item $y \in \ddag{}(V)^{\psi=1}_g$ si et seulement si pour tout $n \geq
%0$, on a $h^1_{K_n,V}(y) \in H^1(K_n,\Fil^1 V)$;
\item si $V$ n'a pas de sous-quotient
stable par $H_K$, alors 
$\ddag{}(V)^{\psi=1}_g = \ddag{}(\Fil^1 V)^{\psi=1}$; 
\item en g\'en\'eral, 
$\ddag{}(\Fil^1 V)^{\psi=1} \subset \ddag{}(V)^{\psi=1}_g$ 
et le quotient est un $\Qp \otimes_{\Zp} \Lambda_K$-module de torsion.
\end{enumerate}
\end{theo}

\begin{proof}
Le fait que $\ddag{}(\Fil^1 V)^{\psi=1} \subset
\ddag{}(V)^{\psi=1}_g$ suit du fait (d\'emontr\'e en 
\cite[lemme 6.5]{LB02} par exemple) que
si $W$ est une repr\'esen\-tation
de de Rham dont tous les poids 
de Hodge-Tate sont $\geq 1$, 
alors toute extension de $\Qp$ par $W$ 
est elle-m\^eme de de Rham. 

Nous allons maintenant montrer le (2), c'est-\`a-dire que 
$\ddag{}(V)^{\psi=1}_g / \ddag{}(\Fil^1 V)^{\psi=1}$ 
est un $\Qp \otimes_{\Zp} \Lambda_K$-module de torsion.
Par un argument de d\'evissage, 
il suffit de montrer que pour tout $k \geq 0$, 
$\ddag{}(\Fil^{-k} V)^{\psi=1}_g /
\ddag{}(\Fil^{1-k} V)^{\psi=1}_g$ est un $\Qp \otimes_{\Zp}
\Lambda_K$-module
de torsion ce que nous allons
maintenant faire.
On a une suite exacte:
\begin{equation}\tag{$S_k$}\label{eqn2}
0 \ra \ddag{}(\Fil^{1-k} V)^{\psi=1}_g \ra \ddag{}(\Fil^{-k} V)^{\psi=1}_g \ra
\ddag{}(\Fil^{-k} V/\Fil^{1-k} V)^{\psi=1}_g 
\end{equation} 
et la repr\'esentation $\Fil^{-k} V/\Fil^{1-k} V$ a des poids $\geq -k$ et
n'admet pas de sous-repr\'e\-sentation dont les poids sont $\geq 1-k$
par le lemme \ref{no_fil_fil}. La
proposition \ref{univ_irred} montre que  
$\ddag{}(\Fil^{-k} V/\Fil^{1-k} V)^{\psi=1}_g \subset 
(\Fil^{-k} V/\Fil^{1-k} V)^{H_K}$ et donc que 
$\ddag{}(\Fil^{-k} V)^{\psi=1}_g  / \ddag{}(\Fil^{1-k} V)^{\psi=1}_g$
est un $\Qp \otimes_{\Zp} \Lambda_K$-module de torsion.
Ceci montre le (2).

Les calculs ci-dessus montrent de plus
que si $V$ n'a pas de sous-quotient
stable par $H_K$, alors en fait $\ddag{}(\Fil^1 V)^{\psi=1} =
\ddag{}(V)^{\psi=1}_g$ car dans ce cas, $\ddag{}(\Fil^{-k}
V/\Fil^{1-k} V)^{\psi=1}_g = 0$ pour tout $k \geq 0$
et donc $\ddag{}(\Fil^{-k} V)^{\psi=1}_g  = \ddag{}(\Fil^{1-k} V)^{\psi=1}_g$
pour tout $k \geq 0$. Ceci montre le (1).

%Montrons le (1).
%La proposition \ref{univ_irred} montre que
%\[ \ddag{}(\Fil^{-k} V/\Fil^{1-k} V)^{\psi=1}_g = \begin{cases}
%(\Fil^{-k} V/\Fil^{1-k} V)^{H_K} & \text{si $k \geq 1$} \\ 
%0 &  \text{si $k = 0$,} \end{cases} \]
%ce qui fait que l'image de $\ddag{}(\Fil^{-k} V/\Fil^{1-k} V)^{\psi=1}_g$
%par $h^1_{K_n}$ est nulle.
%On a donc pour tous $k,n \geq 0$ un morceau de diagramme commutatif:
%\[ \begin{CD}
% \ddag{}(\Fil^{1-k} V)^{\psi=1}_g @>>> \ddag{}(\Fil^{-k}
%V)^{\psi=1}_g @>>> \ddag{}(\Fil^{-k} V/\Fil^{1-k} V)^{\psi=1}_g \\
%@VV{h^1_{K_n}}V  @VV{h^1_{K_n}}V @VV{h^1_{K_n}=0}V \\
%H^1(K_n,\Fil^{1-k} V) @>>> H^1(K_n,\Fil^{-k} V) @>>> H^1(K_n,\Fil^{-k} V/\Fil^{1-k} V)
%\end{CD} \]
%dont on d\'eduit  par une petite chasse au diagramme que
%l'application naturelle
%\[ h^1_{K_n,V}(\ddag{}(\Fil^{1-k} V)^{\psi=1}_g) \ra h^1_{K_n,V}(\ddag{}(\Fil^{-k}
%V)^{\psi=1}_g) \] est surjective, ce qui montre (puisque $V= \Fil^{-k}
%V$ si $k \gg 0$) que l'application
%\[ h^1_{K_n,V}(\ddag{}(\Fil^1 V)^{\psi=1}) \ra h^1_{K_n,V}(\ddag{}(V)^{\psi=1}_g) \]
%est surjective et donc que si $y \in \ddag{}(V)^{\psi=1}_g$ alors pour tout $n \geq
%0$, on a $h^1_{K_n,V}(y) \in H^1(K_n,\Fil^1 V)$. Ceci montre le (1).
\end{proof}

\begin{rema}
Reprenons la d\'emonstration du th\'eor\`eme
\ref{main_in} ci-dessus. L'image de 
$\ddag{}(\Fil^{-k} V)^{\psi=1}_g$ dans $\ddag{}(\Fil^{-k} V/\Fil^{1-k}
V)^{\psi=1}_g$ s'identifie \`a une sous-repr\'esentation $W_{-k}$ de 
$(\Fil^{-k} V/\Fil^{1-k} V)^{H_K}$; c'est 
donc une repr\'esentation de $G_K$ fix\'ee par 
$H_K$ et dont le seul poids de Hodge-Tate
est $-k$ (notons aussi que $W_0=0$). Comme il n'y
a pas d'extensions non-triviales entre de tels objets qui soit encore
fix\'ee par $H_K$ (ou, ce qui revient au m\^eme, il n'y a pas
d'extensions non-triviales entre de tels objets chez les 
$\Qp \otimes_{\Zp} \Lambda_K$-modules), 
on en d\'eduit 
que si l'on {\og splice \fg} les suites exactes
(\ref{eqn2}) pour $k \geq 0$, on trouve:
\[ 0 \ra \ddag{}(\Fil^1 V)^{\psi=1} \ra \ddag{}(V)^{\psi=1}_g \ra
\oplus_{k \geq 1} W_{-k} \ra 0. \]
ce qui permet de pr\'eciser (au moins en principe) le th\'eor\`eme \ref{main_theo}.
\end{rema}

\appendix

\section{Liste des notations}

Voici une liste des principales notations du texte, dans l'ordre o\`u
elles apparaissent:

\ref{I1}: $k$, $K$, $F$, $\mu_{p^n}$, $K_n$, $K_\infty$, $G_K$, $H_K$,
$\chi$, $\Gamma_K$, $F'$, $\sigma$, $\bdag{,r}_F$, $\phi$,
$\bdag{}_F$, $\mathbf{B}_F$, $e_K$, $\bdag{}_K$, $\mathbf{B}_K$,
$\bdag{,r}_K$, $D^\dag$, $\psi$, $r(K)$, $\iota_n$, $D^{\dag,r}$,
$r(D)$. 

\ref{secpgdr}: $\bnrig{,r}{,K}$, $t$, $\bnrig{}{,K}$, $\nabla$, $\partial$,
$\ddagrig$, $\nabla_D$, $\partial_D$, $N_r$, $N$, $\ell_Y$.

\ref{II1}: $V$, $d$, $\ddag{}(V)$, $r(V)$, $\ndr(V)$, $\dcris(V)$,
$\dst(V)$, $\delta_{V(-k)}$.

\ref{II2}: $h^1_{K_n,V}$, $H^1_{\Iw}(K,V)$, $\partial_V$, $r_{\psi}(V)$.

\ref{II3}: $\ddag{}(V)^{\psi=1}_g$, $\Fil^j$.

\end{document}